\documentstyle[11pt,epsfig]{article}
\newcommand{\qed}{\mbox{$\Diamond$}\vspace{\baselineskip}}

\newcommand{\vanish}[1]{}

\begin{document}
\setlength{\baselineskip}{1.2\baselineskip}

\title{An Infinite Antichain of Permutations}
\author{Daniel A. Spielman
        \\
        Department of Mathematics \\
        Massachusetts Institute of Technology \\
        Cambridge, MA 02139
\and
       Mikl\'os B\'ona \\
        Department of Mathematics \\
        Massachusetts Institute of Technology \\
        Cambridge, MA 02139
}

\maketitle
\begin{abstract} We constructively prove that the partially ordered set
of finite permutations ordered by deletion of entries contains an
infinite antichain. \end{abstract}

\section{Introduction} 
When considering a partially ordered set with infinitely many
elements, one should wonder whether it contains 
an infinite antichain (that is, a subset in which each pair of
elements are incomparable). 
It is well known that 
all antichains of $N^k$ (where $(x_1,x_2,\cdots ,x_k)\leq
(y_1,y_2,\cdots ,y_k)$ if and only if $x_i\leq y_i$ for $1\leq i \leq k$
) are finite. (See \cite{LOT}). Another basic result is that all
antichains of the 
partially ordered set of the finite words of a finite alphabet are
finite, where $x<y$ if one can delete some letters from $y$ to get $x$.
(This result is due to Higman and can be found in \cite{lothaire}).

In  this paper we examine this question for the partially ordered set
$P$ of finite permutations with the following $<$ relation:
  if $m$ is less than $n$, and
$p_1$ is a permutation of the set $\{1,2,\cdots ,m\} $ and 
$p_2$ is a permutation of the set $\{1,2,\cdots ,n\} $, then $p_1< p_2$
if and only if we can delete $n-m$ elements from $p_2$ 
so that when we re-name the remaining elements according to
  their rank, we obtain $p_{1}$.
For example,
  $1\:3\:2<2\:4\:5\:3\:1$ as we can delete 4  and 1 
  from the latter to get $2\:5\:3$, which becomes $1\:3\:2$
  after re-naming.
Another way to view this relation is that $p_{1} < p_{2}$
  if there are $n-m$ elements of $p_{2}$ that we can delete
so that the
$i$-th smallest  of the remaining elements preceeds exactly $b_i$
elements, where $b_i$ is the number of elements preceeded by $i$ in
$p_1$. In other words, the $i$-th smallest remaining entry of $p_2$ preceeds
the $j$-th one if and only if $i$ preceeds $j$ in $p_1$.  In short,
$p_1<p_2$ if $p_1$ is ``contained'' in $p_2$, that is, there is a
subsequence in $p_2$  in which any two entries relate to each other as
the corresponding entries in $p_1$.  

We would like to point out that any answer to this question would be
somewhat surprising. If there were no infinite antichains in this
partially ordered set, that would be surprising because, unlike the two
partially ordered sets we mentioned in the first paragraph,  $P$ is
defined over an infinite alphabet and the ``size'' of its elements can
be arbitrarily large.  On the other hand, if there is an infinite
antichain, and we will find one, then it shows that this poset
is more complex in this sense than the poset of graphs ordered by the
operations of edge contraction and vertex deletion. (That
  this poset of
graphs does not contain an infinite antichain is a famous
theorem of Robertson and Seymour~\cite{ros,RS}). This is surprising
too, as graphs are usually much more complex than permutations.

\section{The infinite antichain}
We are going to construct an infinite antichain,
  $\{a_{i}\}$. 
The elements of this
antichain will be very much alike; in fact, they will be identical at the
beginning and at the end. Their middle parts will be very similar,
too. These properties will help ensure that no element is contained
in another one. 

Let $a_1=13,12,10,14,8,11,6,9,4,7 ,3,2,1,5$. 
We view $a_1$ as having three parts: a decreasing sequence
of length three at its beginning, a long alternating
permutation starting with the maximal element of the permutation and
ending with the entry 7 at the fifth position from the right (In this
 alternating part odd entries have only even neighbors and vice cersa.
 Moreover, the odd entries and the even entries form two decreasing 
subsequences so that $2i$ is between $2i+5$ and $2i+3$), and
 a terminating subsequence 3 2 1 5. 

\begin{figure}[htbp]
\hspace*{\fill}{\epsfxsize = 6in \epsfbox{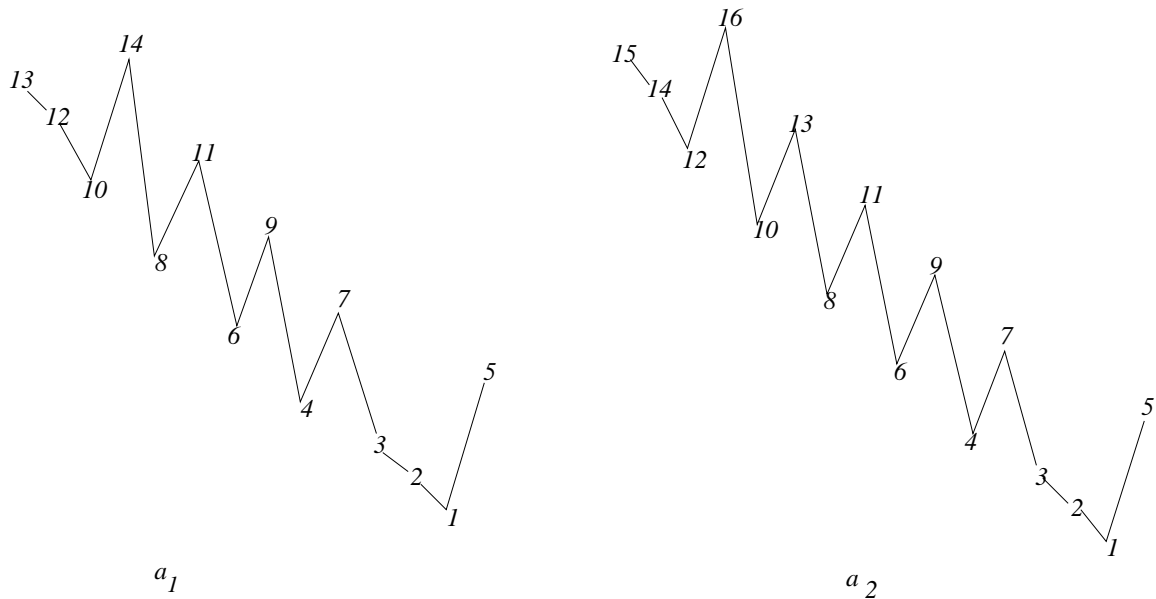}}\hspace*{\fill}
\caption{The pattern of $a_i$}
\label{fig:general}
\end{figure}

To get $a_{i+1}$ from $a_i$, simply insert two consecutive elements
right after the maximum element $m$ of $a_i$, and give them the values
$(m-4)$ and $(m-1)$. Then make the necessary corrections to the rest of
the elements, that is, increment all old entries on the left of $m$ ($m$
included) by two and leave the
rest unchanged (see Figure~\ref{fig:general}). 

Thus the structure of any $a_{i+1}$ is very similar to that of $a_i$---only
the middle part becomes two entries longer.

We claim that the $a_i$ form an infinite antichain. 
Assume by way of contradiction that 
there are indices $i,j$ so that $a_i < a_j$. How could that possibly
happen? First, note that the rightmost element of $a_j$ must map
to the rightmost element of $a_i$, since this is the only element in $a_j$
preceeded by four elements less than itself. 
Similarly, the
maximal element of $a_j$ must map to the maximal element of $a_i$,
since, excluding the rightmost element, this is the only element preceeded by
three smaller elements. This implies that the first four and
the last six elements of $a_j$ must be mapped to the first
  four and last
six elements of $a_i$, thus none of them can be deleted. 

Therefore, when deleting elements of $a_j$ in order to get $a_i$, we can
only delete elements from the middle part, $M_j$. We have already seen
that the maximum element cannot be deleted.  Suppose we can delete a set
$D$ of entries from $M_{j}$ so that the remaining pattern is $a_i$.
First note that $D$ cannot contain three consecutive elements, otherwise
every element before those three elements would be larger than every
element after them, and $a_i$ cannot be divided in two parts with this
property.  Similarly, $D$ cannot contain two consecutive
  elements 
in which the first is even. Thus $D$ can only consist of separate
single elements (elements whose neighbors are not in $D$) and consecutive
pairs in which the first element is odd. 
Clearly, $D$ cannot contain a separate single element as in that case
the middle part of  resulting permutation would contain a decreasing
3-subsequnce, but the middle part, $M_i$, of $a_i$ does not. On the other
hand, if $D$ contained two consecutive elements $x$ and $y$ so that $x$ is
odd, then the odd element $z$ on the right of $y$ would not be in $D$ as
we cannot have three consecutive elements in $D$, therefore $z$ would be
in the remaining copy of $a_i$ and $z$ wouldn't be
preceeded by two entries smaller than itself. This is a contradiction as
all odd entries of $M_i$ have this property. 

This shows that $D$ is necessarily empty, thus we cannot delete any
elements from $a_j$ to obtain some $a_i$ where $i<j$.

We have shown that no two elements in $\{a_i\}$ are comparable,
so $\{a_i\}$ is an infinite antichain. \qed

\end{document}